\documentclass[12pt]{article}
\usepackage[margin=1in]{geometry}
\usepackage[T1]{fontenc}
\usepackage[utf8]{inputenc}
\usepackage{authblk}
\usepackage{amsmath,amssymb}
\usepackage{hyperref}
\usepackage[nameinlink,noabbrev]{cleveref}
\usepackage{graphicx}
\usepackage{caption}

\begin{document}

\title{Statistical warning indicators for abrupt transitions in dynamical systems with slow periodic forcing}

\author[1,*]{Florian Suerhoff}
\author[1,2]{Andreas Morr}
\author[3]{Sebastian Bathiany}
\author[2,3,4]{Niklas Boers}
\author[1]{Christian Kuehn}

\affil[1]{Department of Mathematics, School of Computation, Information and Technology, Technical University of Munich, Boltzmannstraße 3, Garching, Germany}
\affil[2]{Research Domain IV -- Complexity Science, Potsdam Institute for Climate Impact Research, Telegrafenberg A 31, Potsdam, Germany}
\affil[3]{Munich Climate Center and Earth System Modelling Group, Department of Aerospace and Geodesy, TUM School of Engineering and Design, Technical University of Munich, Lise-Meitner-Straße 9, Ottobrunn, Germany}
\affil[4]{Department of Mathematics and Statistics, University of Exeter, North Park Road, Exeter, United Kingdom}
\affil[*]{Correspondence: \href{mailto:f.suerhoff@gmail.com}{f.suerhoff@gmail.com}}

\date{}

\maketitle

\begin{abstract}
There is growing interest in anticipating critical transitions in natural systems, often pursued through statistical detection of early warning signals associated with dynamical bifurcations. In stochastic dynamical systems, such signals commonly rely on manifestations of critical slowing down. However, we still need additional development for the underlying theory for critical transitions in non-autonomous systems. This extension is relevant for natural systems, whose behaviour often emerges from seasonal periodic forcing. In this study, we systematically investigate the feasibility of anticipating the termination of oscillatory behavior in a bistable system with slow periodic forcing. In this setting, existing approaches of estimating linear characteristics of the return map fail in practical scenarios due to the unfavourable time-scale separation. Instead, we propose two statistical indicators for the anticipation of critical transitions in the periodic behaviour: (i) conventional early warning indicators, such as increasing variance and autocorrelation, evaluated across system cycles, and (ii) indicators derived from the phase of the seasonal forcing. By statistically comparing their predictive performance, we find that phase-based indicators provide the strongest early warning capability. Our results offer guidance for the detection of critical transitions in periodically forced systems and, more broadly, systematically extend early-warning signs towards non-autonomous dynamical systems.
\end{abstract}

\section{Introduction}

Abrupt transitions between qualitatively distinct dynamical regimes, often referred to as critical transitions or tipping points, are a central concern across ecology, climate science, and other complex-system settings \cite{ArmstrongMcKay2022TippingPoints, Boers2025DestabilizationTP}. Anticipating such transitions is scientifically and societally relevant because the post-transition regime may be costly or difficult to reverse \cite{Richardson2023PlanetaryBound}. A large body of work studies \emph{early-warning signals} (EWS) based on generic statistical signatures in time series, including rising variance and increasing lag-1 autocorrelation, which are commonly associated with critical slowing down as stability weakens near a bifurcation \cite{scheffer_2009, kuehn_2011, dakos_2012}. Empirical evidence that such indicators can be detectable under field conditions exists in selected settings \cite{carpenter2011_wholeecosystem, Veraart2012RecoveryRatesTipping}, while broader discussions also highlight limitations, confounding effects, and the need for model- and context-aware interpretation \cite{ditlevsen_2010, Morr2024InternalNoiseInterference, Ben-Yami2024TippingTime}.

A key complication in many climate- and ecosystem-relevant applications is that the system of interest is not well described as fluctuating around a static equilibrium. Instead, external drivers impose strong periodic variability (e.g.\ seasonal cycles), and the relevant reference state is a periodic attractor rather than a fixed point \cite{williamson2016_periodicEWS, bathiany2018_oscillatingworld}. This perspective is particularly relevant for prominent Earth-system components that are strongly seasonally forced and have been discussed in terms of nonlinear threshold behavior or abrupt regime shifts, including Arctic sea ice, the North Atlantic subpolar gyre, and major monsoon systems \cite{eisenman2009_seaice, Born2014SPGModel, levermann2009_monsoon}. Abrupt transitions related to quasi-periodic orbital forcing have similarly been discussed under this framework \cite{Imbrie1980ModelingClimatePeriodic, Ditlevsen2009BifurcationStructureNoiseassisted}. Critical transitions in such systems may manifest as a sudden shift in the mean value or even the breakdown of the oscillatory dynamics \cite{Schewe2012MonsoonTransitions, Ben-Yami2024ImpactsAMOCCollapse}.
Accordingly, early-warning methodology that targets stability loss in periodically forced dynamics is directly aligned with the phenomenology emphasized in these application domains \cite{williamson2016_periodicEWS, ma_2018}. In such settings, both the interpretation and the implementation of EWS require adaptation: the statistics of the observed process are strongly phase-dependent within each forcing cycle. This motivates early-warning methodology that is explicitly \emph{cycle-aware}.

In this paper, we study a minimal, analytically tractable model that captures these features: the periodically forced, overdamped Duffing oscillator
\begin{equation}
  \dot x \;=\; x-\tfrac13 x^3 + D_a\cos(\omega t),
  \qquad x(t)\in\mathbb R,
  \label{eq:duffing_det}
\end{equation}
and, when considering stochastic variability, its additive white noise version
\begin{equation}
  dx \;=\; \bigl(x-\tfrac13 x^3 + D_a\cos(\omega t)\bigr)\,dt \;+\; \sigma\,dW_t .
  \label{eq:duffing_sde}
\end{equation}
Here $D_a>0$ controls the forcing amplitude, $\sigma>0$ the noise strength, and $\omega>0$ the forcing rate, with forcing period $T=2\pi/\omega$. The regime of primary interest is slow forcing, meaning a pronounced time-scale separation between the intrinsic (fast) relaxation of the state variable $x$ and the (slow) evolution of the forcing phase. This separation is a natural idealization of seasonally forced systems where within-season adjustment is faster than the seasonal modulation itself \cite{williamson2016_periodicEWS, kuehn2015multiple}.

\begin{figure}
    \centering
    \includegraphics[width=0.8\linewidth]{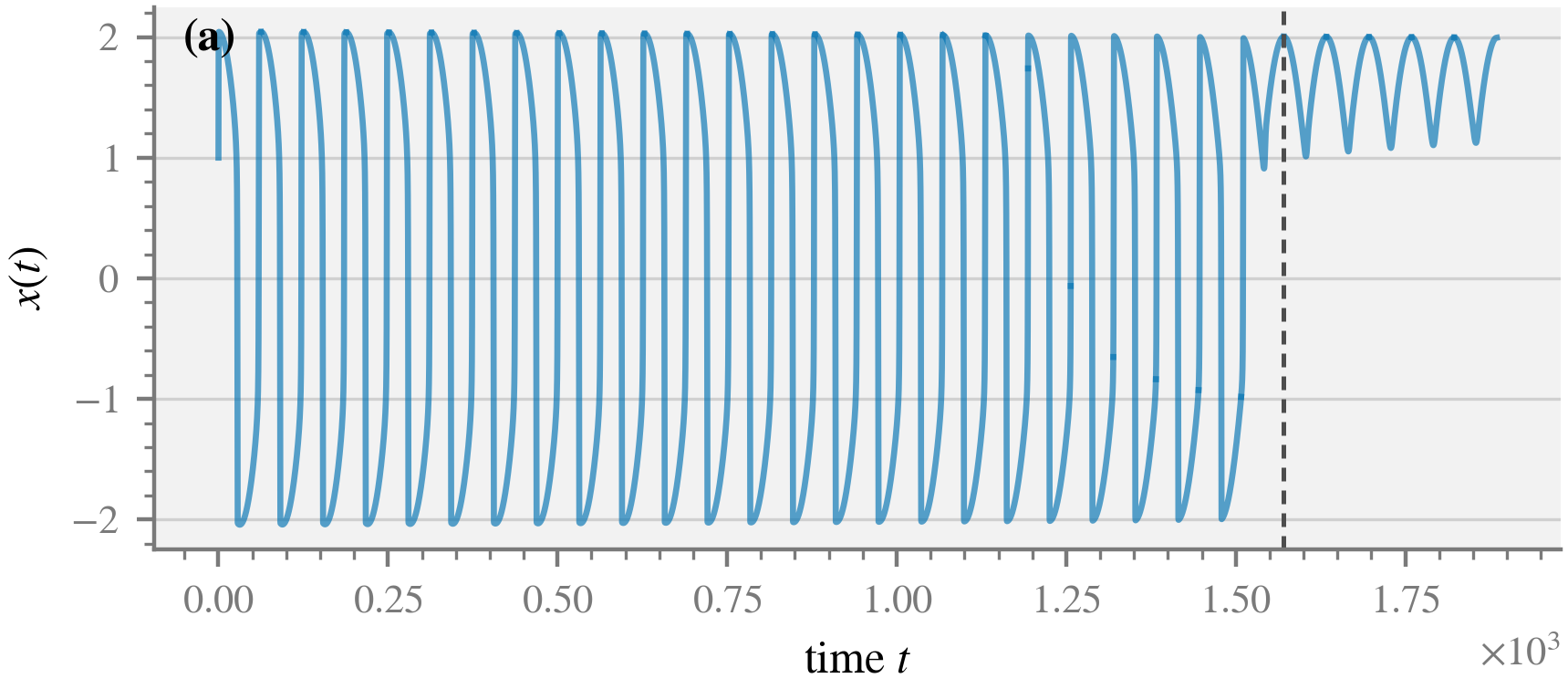}
    \caption{Example realization of the slowly forced Duffing oscillator showing a breakdown of the relaxation oscillation due to a decreasing forcing amplitude $D_a$. The trajectory alternates between wells once per half forcing cycle until a late-cycle transition fails; afterward, the dynamics remain confined to a single-well response with small seasonal modulation. The vertical dashed line marks the breakdown onset (first missed transition).}
    \label{fig:breakdown_timeseries}
\end{figure}

For sufficiently large forcing amplitude, the deterministic dynamics of \eqref{eq:duffing_det} exhibit a fold-mediated relaxation oscillation that jumps between the two wells once per half-cycle. As the amplitude $D_a$ is reduced toward a threshold, this jumping oscillation can break down: the fast jump ceases to occur, and the trajectory remains confined to a single-well response with only small seasonal variability, a regime with undesirable system functioning. This is in contrast with other periodically forced systems investigated in the literature, where the pre-tipping behaviour occurs in a single well, and global transitions to another well occur due to noise-induced phenomena, stochastic resonance, or dynamical bifurcations \cite{Chen2019NTippingPeriodic, berglund_gentz_2002, williamson2016_periodicEWS}. In particular, Williamson et al.~\cite{williamson2016_periodicEWS} seek methods to anticipate such a transition in oscillating systems similar to the one considered here. We here aim to establish such methods for application in systems with repeated relaxation-oscillation jumps.

Figure~\ref{fig:breakdown_timeseries} illustrates the qualitative change under investigation here plotted against time. The central early-warning problem addressed here is to anticipate, under weak stochastic perturbations, this loss of the jumping periodic orbit as $D_a$ approaches the system functioning threshold from above.

\begin{figure}
    \centering
    \includegraphics[width=0.8\linewidth]{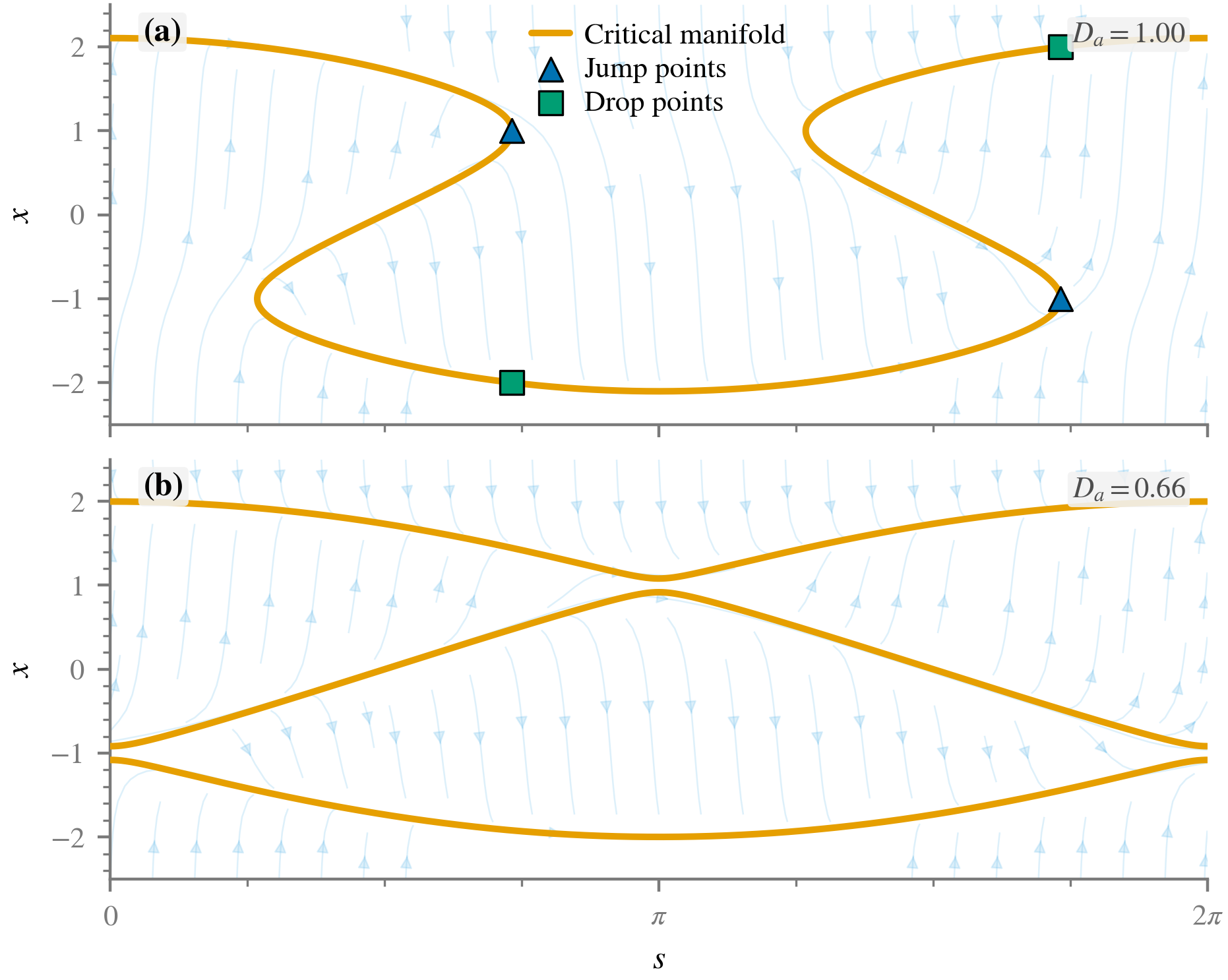}
    \caption{Critical-manifold geometry on the phase cylinder for the fast--slow formulation $\omega\,dx/ds=f(x,s)$ with $f(x,s)=x-\tfrac13x^3+D_a\cos s$. (a) For $D_a>2/3$ the critical manifold has folds and exhibits relaxation oscillations via slow drift along attracting sheets connected by fast jumps; markers indicate representative jump and drop locations. (b) For $D_a<2/3$ the folds vanish and only small-amplitude, single-well seasonal responses persist. Light arrows indicate the reduced flow direction.}
    \label{fig:cylinder_geometry}
\end{figure}

The slow-forcing regime admits a geometric interpretation via a standard fast--slow reformulation. Introducing the phase variable $s=\omega t\in[0,2\pi]$ yields the slow-time form
\begin{equation}
  \omega\,\frac{dx}{ds} \;=\; f(x,s) := x-\tfrac13 x^3 + D_a\cos s,
  \,\omega \ll 1,
  \label{eq:slowfast}
\end{equation}
with associated critical manifold
\begin{equation}
  \mathcal C_0 \;=\; \bigl\{(x,s)\in\mathbb R\times\mathbb S^1:\ f(x,s)=0 \bigr\}.
  \label{eq:critical_manifold}
\end{equation}
The phase-dependent cubic in \eqref{eq:critical_manifold} has folds precisely when $D_a\ge 2/3$, and for $D_a>2/3$ the fast--slow geometry organizes the jumping oscillation into slow passages along attracting sheets connected by fast jumps near the folds \cite{kuehn2015multiple}. The combined cylinder visualization in Figure~\ref{fig:cylinder_geometry} contrasts the folded ($D_a>2/3$) and no-fold ($D_a<2/3$) cases and provides the geometric explanation for the breakdown phenomenon.

A natural first idea for early warning in periodically forced systems is to quantify orbit stability directly from cycle to cycle. For a stable $T$-periodic orbit $x_\omega(t)$, the \emph{Floquet multiplier} $\mu$ (equivalently, the derivative in $x$ of the one-period return map at the periodic orbit) measures linear contraction or expansion over one forcing period \cite{smith_2025}. The multiplier quantifies linear amplification or contraction of small perturbations over one forcing cycle: stability corresponds to $|\mu|<1$, while loss of stability occurs when $|\mu|$ reaches and crosses $1$. In a comparable-time-scale regime (i.e.\ when forcing is not asymptotically slow), one therefore expects a diagnostically useful pattern as the system approaches a transition between dynamical regimes: as the periodic orbit destabilizes (corresponding, in Fig.~\ref{fig:cylinder_geometry}, to the disappearance of the jumping response and a transition from regime~a to regime~b), the cycle-resolved estimate of $\mu$ should increase toward $1$ and then briefly exceed $1$. In other words, in a time series indexed by cycles, $\{\mu_n\}_{n\in\mathbb N}$, the impending regime change would be indicated by a pronounced rise and a spike above $1$ around the critical cycle, making the Floquet multiplier a natural cycle-aware early-warning indicator in that regime. In the slow-forcing limit, however, the relevant relaxation oscillation is typically strongly attracting. In particular, a standard theorem for fast--slow relaxation oscillations (stated for a canonical folded critical-manifold setting) yields linear order dynamics bounded above by $-K/\omega$ and hence a Floquet multiplier bound of the form
\begin{equation}
  |\mu| \;\le\; \exp\!\bigl(-K/\omega\bigr),
  \qquad K>0,
  \label{eq:floquet_bound}
\end{equation}
for sufficiently small $\omega$ \cite{kuehn2015multiple}. While the theorem is formulated under structural assumptions that require minor technical adjustments in the present phase-cylinder geometry, the same fast--slow mechanism applies: as $\omega\to 0$ the orbit spends an overwhelming fraction of each cycle on strongly contracting branches, so the per-period contraction accumulates on the $1/\omega$ time scale and forces $|\mu|$ to become exponentially small. Therefore, the Floquet multiplier scales exponentially and it is in practice extremely difficult to detect a characteristic near-$1$ rise or spike that would make it a practical early-warning diagnostic under slow forcing.

This has an important methodological consequence. In the slow-forcing regime relevant for seasonally forced tipping elements, the Floquet multiplier does not approach the stability threshold $1$ in a gradual, statistically exploitable manner. Instead, $|\mu|$ is typically orders of magnitude below $1$ across the pre-transition regime, so that any cycle-to-cycle variation is easily dominated by stochastic fluctuations or observational noise. Consequently, a Floquet-based early warning is ill-suited as a practical indicator of impending breakdown when time-scale separation is strong, motivating the search for alternative, cycle-aware indicators.

The classical critical-slowing-down picture provides a complementary motivation: as stability weakens, local mean reversion slows down, which can manifest statistically as increases in variance and lag-1 autocorrelation \cite{scheffer_2009,dakos_2012,kuehn_2012}. In periodically forced systems, however, such statistics must be computed in a way that respects phase dependence and the geometry of the oscillation. In this work, we implement and compare known and novel indicators that are expected to occur before critical transitions in our modelling setup \cite{williamson2016_periodicEWS,ma_2018}. We focus on cycle-wise indicators that summarize within-cycle dynamics and track their trends across cycles as $D_a$ drifts toward breakdown. Specifically, we consider (i) cycle-averaged critical-slowing-down indicators computed on detrended within-cycle segments, and (ii) indicators derived from the evolving shape of the relaxation oscillation, in particular the timing of the fast jump relative to the forcing phase. The remainder of the paper details these constructions (Methods), evaluates their predictive performance for breakdown detection under weak noise (Results), and discusses why geometry-informed, phase-based indicators can be especially effective in this slow-forcing setting (Discussion).

\section{Methods}
We construct cycle-wise statistical indicators for critical transitions in the stochastic Duffing oscillator. For this, we consider a single observation time series and extract a small set of scalar summaries from each forcing cycle, whose evolution is tracked across cycles. Therefore, rather than estimating stability from the full return map and its Floquet multipliers, we compress within-cycle information into one value per period to obtain a trendable sequence of indicators.

We first establish the following notation. Throughout this study, we denote the periodic forcing by
\[
\lambda(t)=D_a\cos(\omega t), \qquad T_f=\frac{2\pi}{\omega},
\]
and let $x(t)$ be a (noisy) trajectory of the model. We index all event times by the
\emph{jump number} $k\in\mathbb N$, i.e.\ $k$ enumerates successive well-to-well transitions.

Let $\{t_J^{(k)}\}_{k\ge 1}$ be the strictly increasing sequence of jump times (the fast transitions between wells).
Operationally, we detect jumps using a simple hysteresis rule with two fixed thresholds $x_{\mathrm{up}}>0>x_{\mathrm{low}}$:
starting from a trajectory segment in the upper well, the next down-jump time is defined as the first time $t$ at which
$x(t)$ falls below $x_{\mathrm{low}}$; after such a down-jump, the next up-jump time is the first time $t$ at which
$x(t)$ rises above $x_{\mathrm{up}}$. We then alternate this procedure, yielding an increasing sequence
$t_J^{(1)}<t_J^{(2)}<\cdots$ of detected jump times.
For each jump time $t_J^{(k)}$, we define the associated forcing extremum time
\begin{equation}
t_\star^{(k)} \;:=\; \frac{m_k\pi}{\omega},
\qquad
m_k \in \arg\min_{m\in\mathbb Z}\left|t_J^{(k)}-\frac{m\pi}{\omega}\right|,
\label{eq:textremum_def}
\end{equation}
i.e.\ $t_\star^{(k)}$ is the nearest extremum of $\cos(\omega t)$ (a maximum if $m_k$ is even, a minimum if $m_k$ is odd).
We write $\eta_k:=\cos(\omega t_\star^{(k)})\in\{+1,-1\}$ to indicate whether the associated extremum is a maximum ($\eta_k=+1$) or a minimum ($\eta_k=-1$).

For $D_a>2/3$, define the fold time $t_{\mathrm{fold}}^{(k)}$ as the unique time in the half-cycle leading up to $t_\star^{(k)}$
at which the forcing attains the corresponding static fold value $\pm 2/3$:
\begin{equation}
\lambda\!\left(t_{\mathrm{fold}}^{(k)}\right)=\eta_k\,\frac{2}{3},
\qquad
t_{\mathrm{fold}}^{(k)}\in\Bigl(t_\star^{(k)}-\frac{\pi}{\omega},\,t_\star^{(k)}\Bigr).
\label{eq:tfold_def}
\end{equation}
Equivalently, $\cos(\omega t_{\mathrm{fold}}^{(k)})=\eta_k\,\frac{2}{3D_a}$.

To compare event times across jumps, we map time differences to phase lags by multiplying with $\omega$
and applying the principal-value wrap $\mathrm{wrap}(\cdot):\mathbb R\to(-\pi,\pi]$.
We define the jump phase relative to the associated forcing extremum,
and the dynamic delay relative to the associated fold, by
\[
\Psi^{(k)} := \omega(t_J^{(k)}-t_\star^{(k)}),
\qquad
\Phi^{(k)} := \omega\bigl(t_J^{(k)}-t_{\mathrm{fold}}^{(k)}\bigr).
\]
These satisfy the exact identity
\begin{equation}
\label{eq:decomp}
\Psi^{(k)}
=
\omega\bigl(t_{\mathrm{fold}}^{(k)}-t_\star^{(k)}\bigr)
+
\Phi^{(k)} ,
\end{equation}
where the first term is purely determined from the position of the forcing extremum and the forcing's crossing of the critical threshold. The derivations below are used in a motivational sense: they justify the
phase-based indicator choices and their expected monotone trends as $D_a\downarrow 2/3$.
To support reproducibility and remove ambiguity in how these quantities are computed from simulated trajectories,
Section~\ref{sec:numerical_setup} specifies the numerical protocol and the operational definitions used for jump detection,
segmentation, detrending, and feature construction, including the four final trend features used in the statistical analysis.

\begin{figure}
    \centering
    \includegraphics[width=0.8\linewidth]{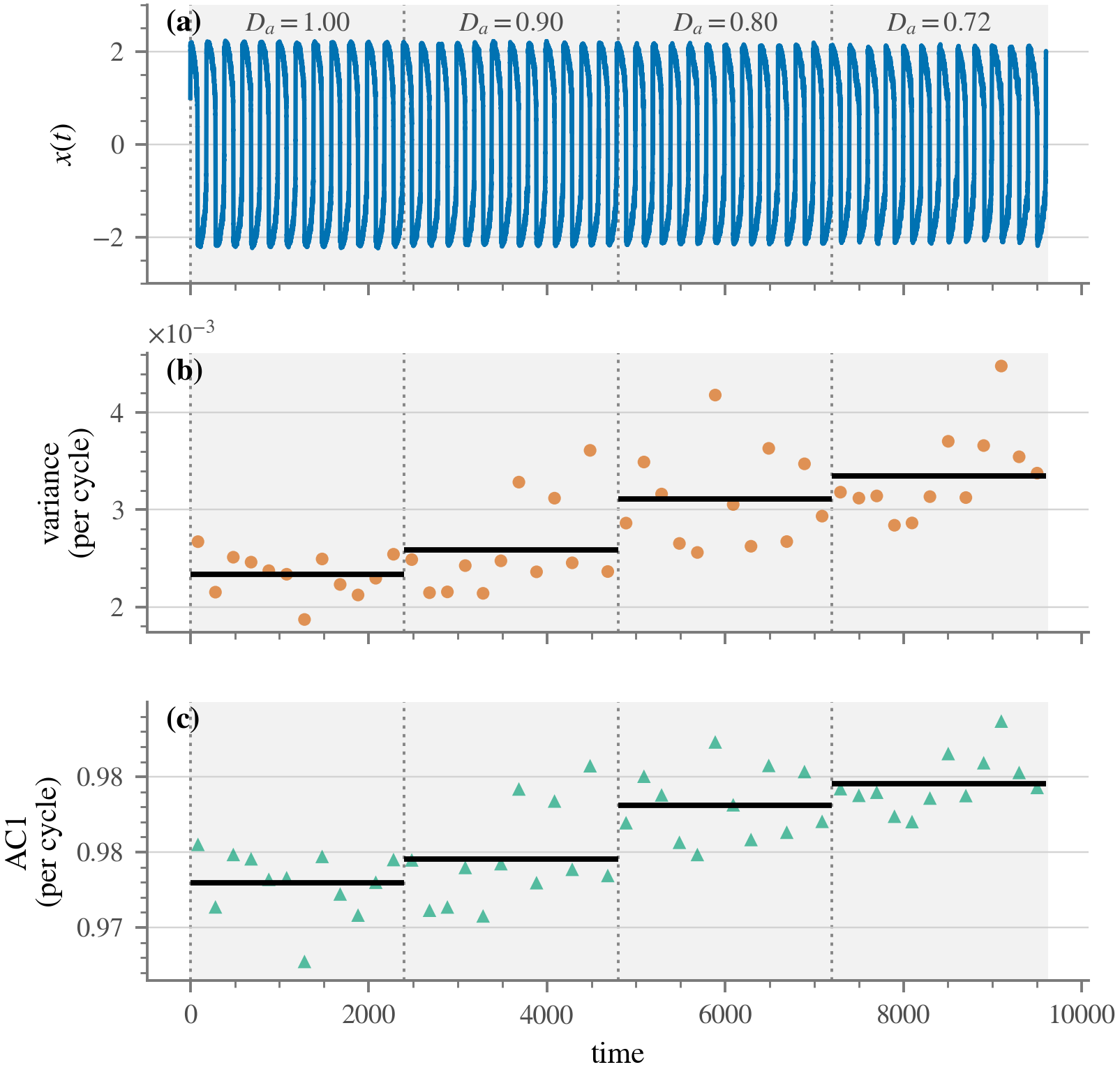}
    \caption{Illustration of cycle-averaged fluctuation indicators under piecewise-constant forcing amplitudes. (a) Example trajectory $x(t)$ for successive amplitude levels $D_a\in\{1.00,0.90,0.80,0.72\}$. (b) Per-cycle variance $\mathrm{Var}_n$ computed from the concatenated detrended residuals of the two between-jump segments in each cycle. (c) Corresponding per-cycle lag-1 autocorrelation $\mathrm{AC1}_n$ computed from the same concatenated residual series. Vertical dotted lines mark amplitude changes; thick black bars denote the mean of the estimated cycle-wise quantities across the time stretch of constant amplitude.}
    \label{fig:var_acf_multilevel}
\end{figure}

\subsection{Cycle-averaged critical slowing down indicators}
Classical early-warning signals for approaching a bifurcation indicate critical slowing down: weaker mean reversion implies increased short-lag memory and increased variance \cite{scheffer_2009,dakos_2012,kuehn_2012}. Often, this is motivated by regarding the nonlinear system as close to equilibrium and employing an Ornstein--Uhlenbeck process to represent the approximately linear dynamics
\[
dy=-\kappa y\,dt+\sigma\,dW_t.
\]
Sampling this process at time steps $\Delta t$ yields an AR(1) process with coefficient $\phi=e^{-\kappa\Delta t}$, so
\[
\mathrm{AC1}=\phi=e^{-\kappa\Delta t}\nearrow 1 \ \text{as }\kappa\downarrow 0,
\qquad
\mathrm{Var}(y)=\frac{\sigma^2}{2\kappa}\nearrow \infty \ \text{as }\kappa\downarrow 0.
\]
This is the standard characteristic system behaviour exploited in a sliding-window routine with prior time series detrending, to ensure the residual represents the mean reversion behaviour. The trend of the statistical quantities is then statistically assessed \cite{kuehn_2011, dakos_2012}.

In the nonautonomous Duffing setting, within-cycle dynamics are strongly phase-dependent and hence nonstationary.
To recover a meaningful notion of local mean reversion, we segment time by jump events and compute fluctuation
statistics on between-jump segments. With $\{t_J^{(k)}\}_{k\ge 1}$ the jump times, define the $k$-th
between-jump segment as
\[
S_k := [t_J^{(k)},\,t_J^{(k+1)}].
\]
On each segment $S_k$ we detrend $x$ and write $y=x-x_{\mathrm{seg}}$, where $x_{\mathrm{seg}}$ is a smooth within-segment trend
(in practice implemented by low-order polynomial regression, consistent with standard EWS preprocessing \cite{dakos_2012}).
Linearizing the drift along the observed trajectory yields the inhomogeneous mean-reverting approximation
\[
\dot y = a(t)\,y + \sigma\,\xi(t),
\qquad
a(t)=\partial_x f\bigl(x(t),\lambda(t)\bigr)=1-x(t)^2,
\]
and we define the segment mean-reversion strength by
\[
\kappa_{S_k} := -\frac{1}{|S_k|}\int_{S_k} a(t)\,dt \;>\;0 .
\]
For a discrete sampling step $\Delta t$ that resolves the segment dynamics, the OU analogy motivates
\[
\mathrm{AC1}[S_k]\approx e^{-\kappa_{S_k}\Delta t},
\qquad
\mathrm{Var}[S_k]\approx \frac{\sigma^2}{2\kappa_{S_k}}.
\]
 As $D_a \searrow 2/3^{+}$, pre-jump portions of the attractor spend longer near $x \approx \pm 1$, where $a(t)=1-x^2 \approx 0$, so the effective mean-reversion strength $\kappa_{S_k}$ decreases on the affected segments. Accordingly, $\mathrm{AC1}[S_k]$ increases toward $1$ and $\mathrm{Var}[S_k]$ increases on those segments. The corresponding increases in the cycle-wise summaries $\mathrm{AC1}_n$ and $\mathrm{Var}_n$, as measured by their linear $\text{slope}_{\mathrm{Var}}$ and $\text{slope}_{\mathrm{AC1}}$ obtained through regression (details below), provide a cycle-resolved, segment-based meta--EWS of weakening restoring forces.

An example trajectory together with the corresponding cycle-averaged $\mathrm{Var}_n$ and $\mathrm{AC1}_n$ trends is shown in Figure~\ref{fig:var_acf_multilevel}.

\begin{figure}
    \centering
    \includegraphics[width=0.8\linewidth]{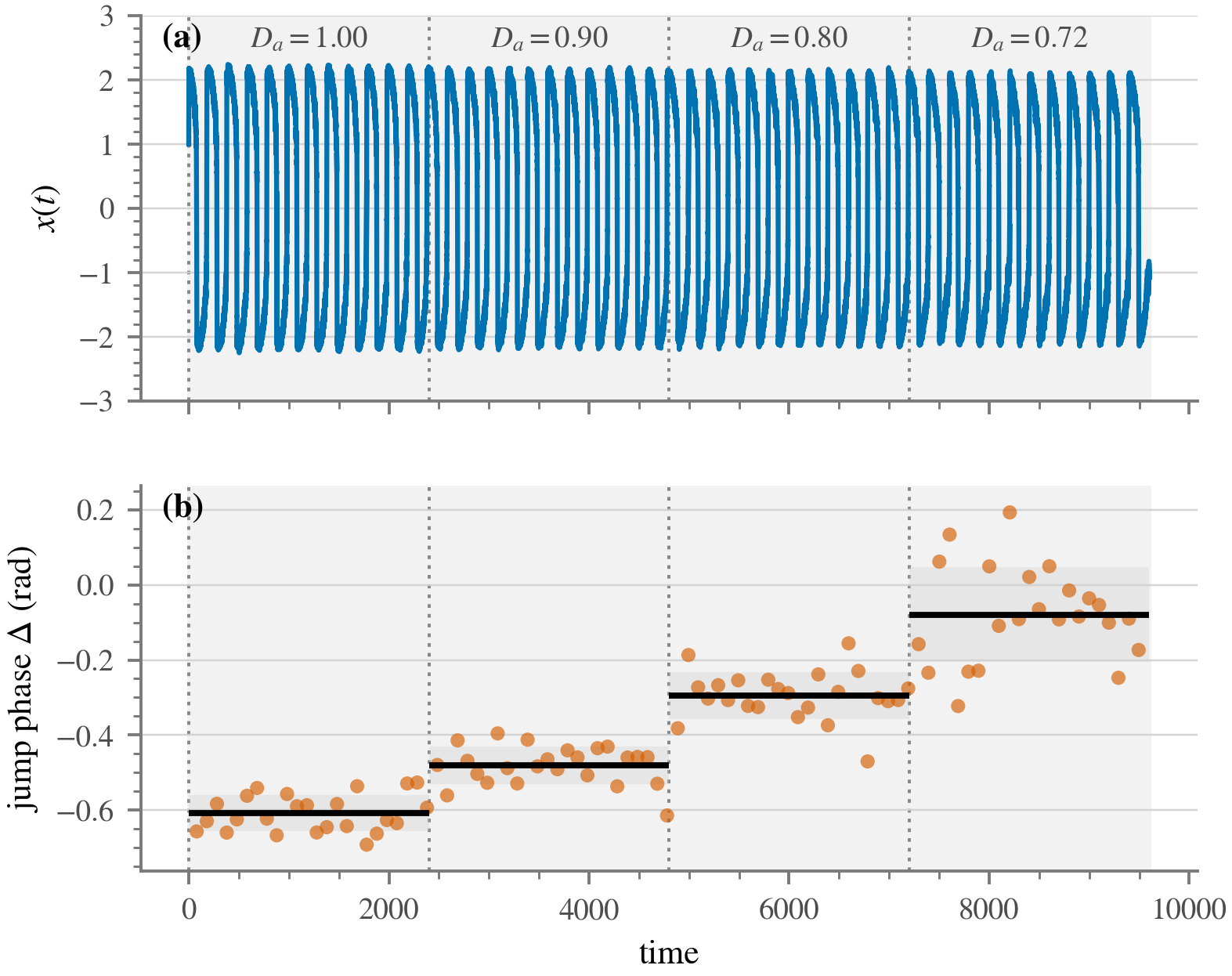}
    \caption{Illustration of phase-based indicators under the same piecewise-constant amplitude protocol as in Figure~\ref{fig:var_acf_multilevel}. (a) Example trajectory $x(t)$. (b) Per-jump phase $\Delta_j$ (signed phase difference to the nearest forcing extremum) together with within-level circular mean (black) and dispersion (shaded band). As $D_a$ decreases toward the fold threshold, the jump phase drifts toward the associated extremum and the jump-phase spread increases.}
    \label{fig:mean_std_multilevel}
\end{figure}

\subsection{Cycle-shape based indicators}
The phase-based indicators exploit the decomposition \eqref{eq:decomp}, which splits the jump phase $\Psi^{(k)}$
(relative to the associated forcing extremum $t_\star^{(k)}$) into an offset determined by the forcing's threshold crossing and its maximum,
and a dynamical delay $\Phi^{(k)}$ induced by slow passage near the fold.
By definition \eqref{eq:tfold_def}, the fold time $t_{\mathrm{fold}}^{(k)}$ lies one half-cycle before $t_\star^{(k)}$ and satisfies
$\cos(\omega t_{\mathrm{fold}}^{(k)})=\eta_k\,\frac{2}{3D_a}$ with $\eta_k=\cos(\omega t_\star^{(k)})\in\{\pm 1\}$.
Hence, the phase distance between the fold and the associated extremum is the same for maxima and minima, giving the first term of \eqref{eq:decomp} as
\[
\omega\bigl(t_{\mathrm{fold}}^{(k)}-t_\star^{(k)}\bigr)
=-\arccos\!\Bigl(\frac{2}{3D_a}\Bigr)\in(-\pi,0),
\]
which increases to $0$ as $D_a\downarrow 2/3$ (the fold moves toward the associated forcing extremum, consistent with the critical-manifold geometry shown in Fig.~\ref{fig:cylinder_geometry}).

For the dynamical delay $\Phi^{(k)}$, consider the local saddle--node (fold) normal form viewpoint and define the fold-distance parameter 
\[
\delta_{\mathrm{fold}}(t):=\lambda(t)-\frac23.
\]
Near $t_{\mathrm{fold}}^{(k)}$, one has $\delta_{\mathrm{fold}}(t)\approx \delta_{\mathrm{fold}}'(t_{\mathrm{fold}}^{(k)})(t-t_{\mathrm{fold}}^{(k)})$ with rate
\[
\beta:=\bigl|\delta'(t_{\mathrm{fold}}^{(k)})\bigr|
=\bigl|\lambda'(t_{\mathrm{fold}}^{(k)})\bigr|
= D_a\,\omega\,\sqrt{1-\frac{4}{9D_a^2}}.
\]
A standard slow-passage law for time-dependent saddle--node bifurcations yields
\[
t_J^{(k)}-t_{\mathrm{fold}}^{(k)} = K\,\beta^{-1/3}\qquad(\beta\downarrow 0),
\]
for some $K>0$ \cite{li_2019}. Substituting into $\Phi^{(k)}=\omega\bigl(t_J^{(k)}-t_{\mathrm{fold}}^{(k)}\bigr)$ gives the explicit dependence
\begin{equation}
\Phi^{(k)} = C\,\omega^{2/3}\Bigl(D_a\sqrt{1-\frac{4}{9D_a^2}}\Bigr)^{-1/3},
\qquad C>0,
\label{eq:Phi_scaling}
\end{equation}
showing that $\Phi^{(k)}$ increases as $D_a\searrow 2/3$ (since $\beta\searrow 0$). Together with $\Theta(D_a)\nearrow 0$, the decomposition \eqref{eq:decomp} implies a monotonic increase of $\Psi^{(k)}$ over consecutive jumps as the fold threshold is approached. The jump phase relative to the associated forcing extremum therefore provides a geometry-informed indicator of critical transitions complementary to the fluctuation-based indicators.

Because phases are quantities on a circular axis, we summarize jump timing using a phase representation and circular dispersion.
For each detected jump at time $t_J^{(k)}$, let
\[
\phi_k := (\omega\, t_J^{(k)}) \bmod 2\pi \in [0,2\pi)
\]
be its forcing phase, and assign it to the nearest forcing-extremum phase $\phi^\star\in\{0,\pi\}$ (maximum/minimum of $\cos$).
We then define the signed minimal phase difference
\begin{equation}
    \Delta_k = \mathrm{wrap}\bigl(\phi_k-\phi^\star\bigr)\in(-\pi,\pi],
\end{equation}
where $\mathrm{wrap}:\mathbb{R}\to(-\pi,\pi]$ denotes the principal-value map.

Our first phase-based trend feature is the linear slope of the sequence $\{\Delta_k\}$ versus the jump index $k$; in the numerical analysis below, this feature is denoted by $\text{slope}_{\mathrm{jump\ phase}}$.
To quantify phase dispersion, we compute a rolling circular standard deviation on windows of $W$ consecutive jumps.
For a window of phases $\{\Delta_j\}_{j=1}^n$, define the mean resultant length
\[
R \;=\; \left|\frac{1}{n}\sum_{j=1}^n e^{i\Delta_j}\right|.
\]
The circular standard deviation is then
\begin{equation}
\mathrm{circstd}
=\sqrt{-2\log R}.
\label{eq:circstats}
\end{equation}
The resulting rolling series of circular standard deviations is then summarized by its linear slope versus jump index; in the numerical analysis below, this feature is denoted by $\text{slope}_{\mathrm{phase\ std}}$. Together, these two quantities yield cycle-shape based meta-EWS for systematic drift and increasing timing variability near breakdown.

To motivate the increase in jump-phase dispersion near breakdown, let $\Delta U(t)$ denote the instantaneous potential barrier height between the attracting well and its saddle. Near a fold, $\Delta U(t)\asymp c\,|\delta_{\mathrm{fold}}(t)|^{3/2}$ \cite{berglund_gentz_2002}, and the corresponding instantaneous escape probability obeys a Kramers law
\[
r_\sigma(t)\sim \exp\!\left(-\frac{\Delta U(t)}{\sigma^2}\right)
\qquad\text{for small }\sigma>0
\]
\cite{haenggi_1990}. Using $\delta_{\mathrm{fold}}(t)\approx \pm \beta |t-t_{\mathrm{fold}}^{(k)}|$ near the fold yields
\[
\Delta U(t)\asymp c\bigl(\beta|t-t_{\mathrm{fold}}^{(k)}|\bigr)^{3/2}.
\]
Hence, the time interval in which $r_\sigma(t)$ is $\mathcal O(1)$ is determined by $\Delta U(t)\lesssim C\sigma^2$, equivalently
\[
|t-t_{\mathrm{fold}}^{(k)}|\lesssim C'\,\frac{\sigma^{4/3}}{\beta}.
\]
With the decrease in $D_a\downarrow 2/3$ comes a decrease in  $\beta\downarrow 0$, so this high-hazard window widens like $\mathcal O(\sigma^{4/3}\beta^{-1})$, which manifests as an increasing circular standard deviation of jump phases. In contrast, when $D_a$ is well above the threshold and one well disappears robustly each cycle, slow-passage dynamics lock the jump time close to its deterministic value and the phase dispersion remains approximately constant under weak noise \cite{berglund_gentz_2001}. Figure~\ref{fig:mean_std_multilevel} illustrates the drift of the jump phase and the increase of jump-phase dispersion across decreasing $D_a$ levels.

\section{Results}
\subsection{Experimental setup}
\label{sec:numerical_setup}

We simulate the one-dimensional stochastically forced Duffing-type dynamics
\begin{equation}
    \mathrm{d}x(t) = f\!\bigl(x(t),t;D_a(t)\bigr)\,\mathrm{d}t + \sigma\,\mathrm{d}W_t,
    \qquad
    f(x,t;D_a)=x-\tfrac13 x^3 + D_a \cos(\omega t),
\end{equation}
using an Euler--Maruyama scheme on a uniform grid $t_n=n\Delta t$ with $\Delta t=0.01$ over a total horizon $T_{\mathrm{tot}}=2500$. We use the \emph{angular} forcing rate $\omega=2\pi/225$, so that the forcing period is $T_f=2\pi/\omega=225$. The forcing amplitude is ramped linearly,
\begin{equation}
    D_a(t)=D_{\max}-\bigl(D_{\max}-D_{\min}\bigr)\frac{t}{T_{\mathrm{tot}}},
\end{equation}
with fixed $D_{\max}=1.2$ and realization-specific $D_{\min}\sim \mathrm{Unif}[0.25,0.9]$. Unless stated otherwise, the noise strength is $\sigma=0.3$ and the initial condition is $x(0)=1.0$. We generate $N_{\mathrm{sim}}=1000$ independent realizations.

To robustly identify well-to-well transitions in noisy trajectories, we employ a jump detector with hysteresis thresholds. A discrete state variable tracks whether the trajectory is in the upper or lower well; a jump is recorded when the signal crosses the respective threshold.  Transitions to the lower well have threshold $x_{\mathrm{up}}=0.4$ and transitions to the upper well $x_{\mathrm{low}}=-0.4$. Let $\{i_j\}$ denote the resulting jump indices on the simulation grid, and define the corresponding jump times by $t_J^{(j)}:=t_{i_j}$. We augment the index list with endpoints ($i_0=0$ and $i_{J}=N$) and define between-jump segments as index intervals $(i_j,i_{j+1})$. To avoid micro-segments caused by threshold chatter, we discard segments shorter than $N_{\min}=80$ time steps.

We define a breakdown onset as the first missed transition event, operationalized by an anomalously long between-jump segment. Let $T_f$ be the forcing period; the half-period is $T_f/2$. A realization is labeled as a breakdown if there exists a segment whose duration exceeds $3T_f/4$.
The breakdown onset is the first such segment in time order. For all feature computations, we truncate the data at this onset and only use segments strictly prior to it (non-breakdown runs use all available segments).

All fluctuation-based indicators are computed on detrended within-segment residuals. Given a segment with $n$ samples, we remove a symmetric buffer fraction $p_{\mathrm{buf}}=0.05$ of the time series on both ends to reduce contamination from the fast transition dynamics near jumps. On the remaining interior samples, we fit a cubic polynomial trend and define residuals as the detrended signal.

We group consecutive between-jump segments into forcing cycles by pairing $(\text{seg}_{2n-1},\text{seg}_{2n})$,
where $n$ denotes the forcing-cycle index. For each cycle, we concatenate the residuals from the two paired segments
and compute (i) the sample variance $\mathrm{Var}_n$, and (ii) the lag-one autocorrelation $\mathrm{AC1}_n$. Trend features are defined as ordinary-least-squares (OLS) slopes of $\{\mathrm{Var}_n\}$ and
$\{\mathrm{AC1}_n\}$ versus the cycle index $n$, denoted by $\text{slope}_\mathrm{Var}$ and $\text{slope}_\mathrm{AC1}$. Slopes are only estimated when at least five cycles are available; otherwise
the run is excluded from the downstream analysis.

For each detected jump time $t_J^{(j)}$ (excluding the artificial endpoints), we compute the forcing phase
\begin{equation}
    \phi_j = (\omega\, t_J^{(j)}) \bmod 2\pi \in [0,2\pi).
\end{equation}
We assign the jump to the nearest forcing extremum phase $\phi^\star\in\{0,\pi\}$ (corresponding to a maximum or minimum of $\cos$) and define the signed minimal phase difference
\begin{equation}
    \Delta_j = \mathrm{wrap}\bigl(\phi_j-\phi^\star\bigr)\in(-\pi,\pi],
\end{equation}
where $\mathrm{wrap}(\cdot)$ maps angles to $(-\pi,\pi]$. The first phase trend feature is the OLS slope of $\{\Delta_j\}$ versus the jump index $j$, denoted by $\text{slope}_{\mathrm{jump\ phase}}$ (computed only if at least five jumps are available). To quantify drift in jump-phase dispersion, we compute a rolling standard deviation  $s_j$ (computed on the circular axis) from $\{\Delta_j\}$ using a window of $W=16$ consecutive jumps. The second phase trend feature is the OLS slope of $\{s_j\}$ versus the jump index $j$, denoted by $\text{slope}_{\mathrm{phase\ std}}$.

Each realization thus yields four trend features: the variance slope, the $\mathrm{AC1}$ slope, the jump-phase slope, and the jump-dispersion slope; the breakdown label is defined by the criterion above.

For the supervised benchmark, we use a linear SVM trained on standardized features via a pipeline (feature standardization performed within each cross-validation split). Concretely, for each fold, we fit the scaler on the training portion only, transform both train and validation features with that fitted scaler, and then fit the linear SVM on the transformed training data to prevent any information leakage from the validation set. We report balanced accuracy using 5-fold stratified cross-validation with fixed random seeds for the split and all stochastic components.

\subsection{Quality of the statistical indicators}

We assess whether the four cycle-level trend features
\[
\bigl(\text{slope}_\mathrm{Var},\ \text{slope}_{\mathrm{AC1}},\ \text{slope}_{\mathrm{mean\ phase}},\ \text{slope}_{\mathrm{phase\ std}}\bigr)
\]
carry predictive information about an oscillation breakdown of the relaxation oscillation under a global amplitude ramp. We fix the forcing frequency in the slow-forcing regime and simulate \(R=1000\) independent realizations of the noisy Duffing oscillator. In each run, the forcing amplitude follows a linear ramp form, with a fixed start value $D_{\max}$ and an independently drawn final value $D_{\min}$ in each run, yielding a mixture of trajectories that keep cycling regularly over the full horizon and trajectories that lose the alternating-cycle structure close to the end of the ramp. This assesses the indicators' capabilities in varying settings, a faster and more substantial forcing change when $D_{\min}$ is small, and a gradual and less impactful forcing change when it is large. Each trajectory is simulated over approximately eleven periods of the periodic driver.

Each realization is assigned the label "breakdown" or "no breakdown" by the segment-duration criterion: the first between-jump segment exceeding \(1.5\) half-periods marks the onset of breakdown, and data after this onset are excluded from feature construction. For every run, the four indicators are computed as slopes across forcing cycles, with variance and \(\mathrm{AC1}\) extracted from detrended, buffered residual segments \cite{scheffer_2009,dakos_2012}, and phase statistics extracted from the detected jump times.

\begin{figure}[t]
  \centering
  \includegraphics[width=0.8\linewidth]{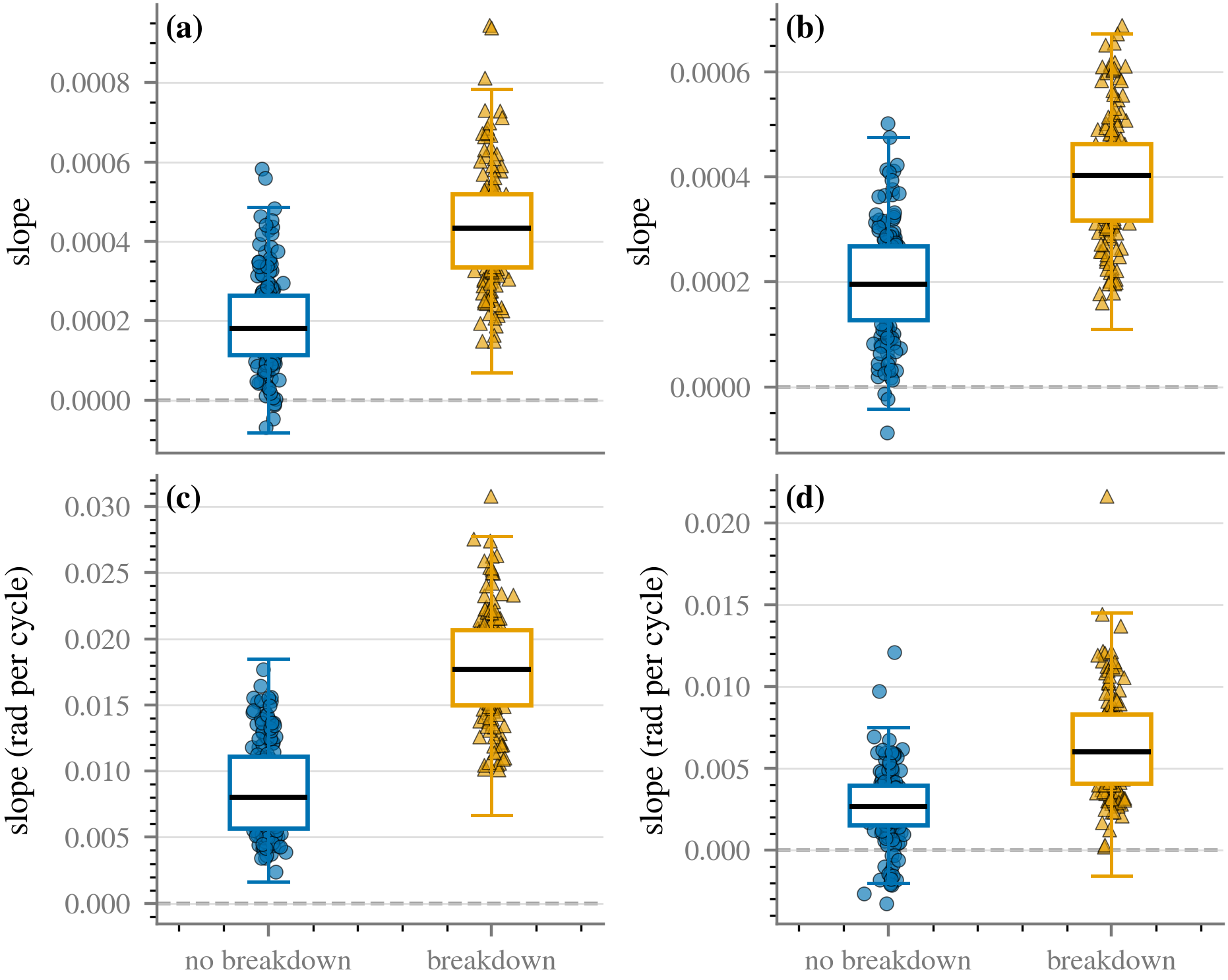}
  \caption{Indicator slopes grouped by whether the respective trajectory exhibited a breakdown, shown for a subsample of 150 runs. The indicators are (a) \(\text{slope}_\mathrm{Var}\), (b) \(\text{slope}_{\mathrm{AC1}}\), (c) \(\text{slope}_{\mathrm{mean\ phase}}\), (d) \(\text{slope}_{\mathrm{phase\ std}}\). The mean-phase slope provides the clearest class separation, followed by \(\mathrm{AC1}\); variance and phase dispersion exhibit stronger overlap.}
  \label{fig:boxplots-byclass}
\end{figure}

Figure~\ref{fig:boxplots-byclass} summarizes the class-conditional distributions of these per-run slopes (visualized on a stratified subsample of 150 runs for readability). The segment-based indicators \(\text{slope}_\mathrm{Var}\) and \(\text{slope}_{\mathrm{AC1}}\) tend to shift upward for breakdown runs, consistent with progressively weakening mean reversion within cycles as the system is driven toward the loss-of-stability region. The separation is visibly stronger for \(\text{slope}_{\mathrm{AC1}}\) than for \(\text{slope}_\mathrm{Var}\), with more overlap in the variance slopes. The most pronounced class separation is observed for the phase-based indicator \(\text{slope}_{\mathrm{mean\ phase}}\): breakdown runs exhibit substantially larger drifts of the mean jump phase across cycles. In contrast, \(\text{slope}_{\mathrm{phase\ std}}\) shows only a weak shift with considerable dispersion.

To quantify joint predictive performance, we train a linear support vector machine (SVM) on the standardized four-feature vectors using stratified 5-fold cross-validation, reporting balanced accuracy to account for class imbalance. The resulting cross-validated balanced accuracy for the full feature set is \(0.873\), indicating that the combined indicator trends provide a strong discriminative signal for breakdown versus no breakdown under the global ramp.

\begin{figure}[t]
  \centering
  \includegraphics[width=0.8\linewidth]{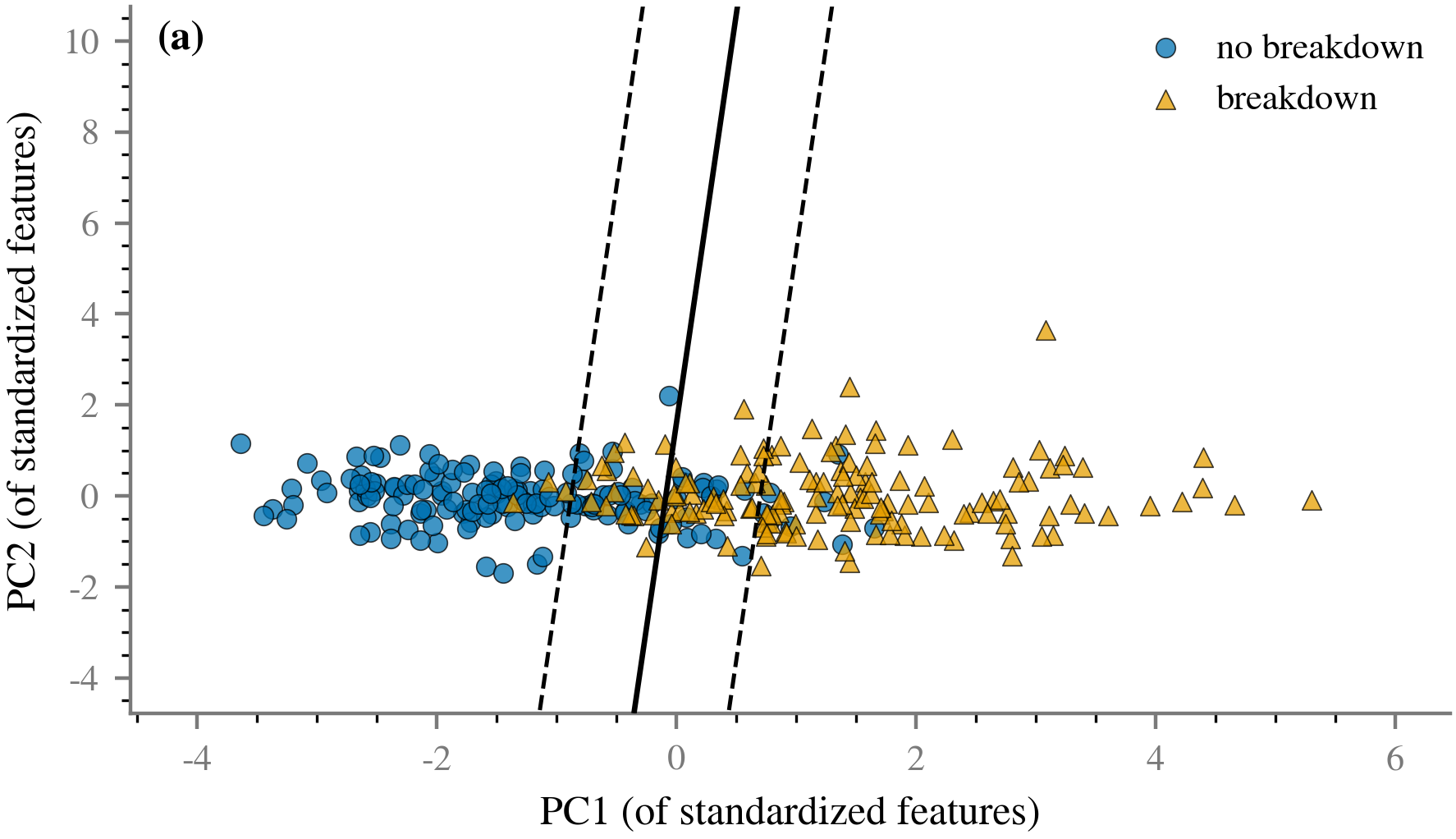}
  \caption{Linear SVM trained on the four indicator slopes; decision geometry projected onto the first two principal components (PC1--PC2) of the standardized feature space. For readability, only 150 runs are plotted; the model is trained and evaluated (stratified 5-fold CV) on all \(R=1000\) runs, achieving balanced accuracy \(0.873\).}
  \label{fig:svm-pca}
\end{figure}

Figure~\ref{fig:svm-pca} visualizes the learned linear separator by projecting the standardized feature space onto its first two principal components; the plot overlays a stratified subsample of 150 runs for readability, while training and scoring use the full \(R=1000\) runs. Consistent with the univariate separations in Figure~\ref{fig:boxplots-byclass}, feature-importance analyses (drop-column \(\Delta\)CV, and permutation importance) identify \(\text{slope}_{\mathrm{mean\ phase}}\) as the most influential predictor, followed by \(\text{slope}_{\mathrm{AC1}}\). In contrast, \(\text{slope}_\mathrm{Var}\) and \(\text{slope}_{\mathrm{phase\ std}}\) contribute only marginally in the presence of the first two features. These importance statements should be interpreted qualitatively, since the four features are partially correlated and multicollinearity can redistribute weight among predictors without materially changing the separator.

\section{Discussion}

This work addresses the anticipation of abrupt changes in periodically forced systems in the slow-forcing regime, where the relevant reference state is a periodic attractor rather than an equilibrium \cite{williamson2016_periodicEWS, bathiany2018_oscillatingworld}. Classical early-warning signals (EWS) based on critical slowing down (notably increasing variance and lag-1 autocorrelation under windowing and detrending) are well-motivated in quasi-static settings \cite{scheffer_2009,dakos_2012,lenton_2012,dakos_2024}, but their interpretation becomes more subtle when statistics vary strongly over the forcing phase \cite{ma_2018,williamson2016_periodicEWS}. In addition, stability metrics derived from the deterministic return map can become structurally uninformative under pronounced time-scale separation: when recovery within each season is fast compared to the forcing period, period-averaged contraction dominates the multiplier and tends to hide the impending loss of the jumping cycle.

In our setting, we introduced cycle-aware indicators that compress within-cycle information into one scalar per cycle and then quantify precritical drift through simple trend features. The indicators are (i) segment-based fluctuation measures (variance and $\mathrm{AC1}$ computed on detrended between-jump segments) and (ii) phase-based measures that explicitly use the timing of fast transitions relative to the forcing phase. In the numerical experiments, these indicators yield high predictive performance for distinguishing runs that exhibit a cycle breakdown from those that maintain regular jumping. In particular, the classification score of $0.873$ is notable given that there is no sharp separation between breakdown and non-breakdown trajectories in the simulated dataset: because the final value $D_{\mathrm{min}}$ is drawn from a uniform distribution, the boundary between the two classes is diffuse, and some runs are labeled as non-breakdown by the segment-duration criterion even though a breakdown appears imminent but would occur only after the end of the simulation window. In a more favourable assessment scenario where $D_{\mathrm{min}}$ is drawn from two well-separated values (one breakdown-inducing and one stable), the class separation becomes correspondingly cleaner. The supervised benchmark further shows that phase-based information improves classification and that the slope of the mean jump phase is the dominant feature among the four trend descriptors. In the drop-column analysis, omitting this feature caused the largest reduction in balanced accuracy (\(\Delta = 4.0\cdot10^{-2}\)), whereas the decreases from omitting the \(\mathrm{AC1}\) slope, the jump-phase standard-deviation slope, and the variance slope were much smaller (\(6.0\cdot10^{-3}\), \(5.0\cdot10^{-3}\), and \(3.0\cdot10^{-3}\), respectively). Permutation importance confirms the same ordering: the mean jump-phase slope attains \((3.2 \pm 0.3)\cdot 10^{-1}\), far exceeding the \(\mathrm{AC1}\) slope at \((4.0 \pm 1.7)\cdot 10^{-2}\), the variance slope at \((7.3 \pm 8.5)\cdot 10^{-3}\), and the jump-phase standard-deviation slope at \((3.3 \pm 7.2)\cdot 10^{-3}\). This indicates that the \(\mathrm{AC1}\) slope supplies only limited complementary information, while variance and jump-phase dispersion make comparatively weak additional contributions once the mean jump-phase trend is present.

A central finding is that indicator strength is linked to how directly the statistic encodes the mechanism that produces the transition. Segment-based variance and $\mathrm{AC1}$ are comparatively agnostic summaries of stochastic fluctuations. In fast--slow systems near a fold, weakened local restoring rates on attracting slow manifolds lead to increasing variance and autocorrelation in appropriately processed residuals \cite{scheffer_2009, kuehn_2011, kuehn_2012, dakos_2012}. In periodically forced dynamics, however, this weakening is phase-dependent and localized; consequently, cycle-averaged fluctuation trends can be diluted by within-cycle heterogeneity and by strong mean reversion away from the near-fold region.

By contrast, the mean jump phase is informative for a more mechanistic reason: it is tied to the geometry of the fold and to the slow-passage dynamics that govern the timing of the fast jump. As the forcing amplitude approaches the fold threshold, the phase at which the potential loses a well shifts toward the forcing maximum, producing a deterministic drift in jump timing. Moreover, dynamic bifurcation delay depends on the local sweep rate through the fold; this rate decreases as the threshold is approached, yielding an additional systematic contribution to delayed escape times \cite{kuehn2015multiple, li_2019}. Under weak noise, stochastic variability broadens the distribution of jump times near the fold, consistent with pathwise metastability and the widening of a high-hazard window for noise-induced tipping as the barrier collapses \cite{berglund_gentz_2001,berglund_gentz_2002,ma_2018}. Importantly, the mean jump phase combines these effects in a way that retains a substantial deterministic component inherited from the slowly varying phase geometry. This reduces the effective noise level of the statistic relative to purely fluctuation-based indicators, which helps explain why mean-phase trends separate breakdown from non-breakdown most cleanly and emerge as the most influential feature in the multivariate classifier.

This observation also makes explicit a trade-off between the needed system knowledge and the obtained indicator ability. Phase-based indicators require (a) access to, or reconstruction of, a meaningful phase variable for the periodic driver, and (b) a reliable definition of the within-cycle event whose timing is tracked (here, the fast jump between wells). When such a structure is available, and the transition mechanism is indeed a fold-mediated jumping, phase-based indicators can be substantially more informative than generic fluctuation summaries. Conversely, when the dynamical mechanism is uncertain, for instance, if transitions are predominantly noise-induced rather than bifurcation-driven, or if rate effects dominate, the predictive value and interpretability of phase-based indicators can degrade, and any apparent warning should be interpreted with caution \cite{ditlevsen_2010, Ashwin_2012}. In that sense, phase-aware indicators are powerful but also more model-committed: their success depends on confidence that the chosen dynamical abstraction captures the relevant geometry and event structure.

Overall, the present results provide a first, systematic assessment of several feasible methods to anticipate critical transitions in slowly forced, periodically driven systems. The main implication is practical: in regimes where period-averaged stability metrics are structurally disadvantaged by time-scale separation, early-warning methodology should shift from per-period return-map characteristics to cycle-wise summaries that respect the structure of the relaxation oscillation \cite{williamson2016_periodicEWS,ma_2018}. The combination of a geometric, phase-based indicator (mean jump phase) with a complementary fluctuation indicator ($\mathrm{AC1}$) appears particularly effective in the studied setting, as it couples mechanistic information about fold timing with stochastic information about weakening mean reversion on attracting slow manifolds.

At the same time, several limitations must be kept in view. First, our analysis is carried out in a minimal, analytically tractable model with additive noise and a specific ramp protocol. While this choice is advantageous for interpretability and analytical footing, empirical systems may exhibit nonstationary or state-dependent noise \cite{Morr2024KramersMoyalEWS}, alternative noise forcing \cite{Riechers2025DiscontinuousDOStatistics, Morr2025RedNoise}, or additional state variables that make the system dynamics more complex \cite{Morr2024InternalNoiseInterference}. Second, the phase-based indicators rely on identifying jumps and on the availability of a meaningful seasonal phase reference; both steps can, in practice, be confounded by measurement noise \cite{BenYami2023AMOCDataCSD}, irregular sampling \cite{Liu2025LambdaLambda}, or ambiguous event definitions \cite{Slattery2024LeadLag}.

Future research should therefore test the proposed indicators under increasingly realistic observation constraints and, ultimately, in real-world periodically forced systems where abrupt breakdown of seasonal cycling is hypothesized \cite{eisenman2009_seaice, levermann2009_monsoon}. Promising directions include: assessing robustness to common data issues (missing cycles, irregular sampling, short records, measurement noise), deriving practical guidelines for windowing and detrending choices, and examining sensitivity to the functional form of the periodic forcing while keeping comparable amplitude ranges. Such validation efforts would help determine to what extent the strong performance of phase-based indicators observed here transfers to complex, seasonally driven systems of applied interest \cite{williamson2016_periodicEWS, bathiany2018_oscillatingworld}. On the modeling side, it would be useful to extend the structural arguments behind the failure of period-aggregated Floquet multipliers under slow forcing beyond the present double-well setting, clarifying when time-scale separation generically renders return-map stability measures poor early-warning statistics and when structurally informed indicators become indispensable.

%
%


\vspace{1em}
\noindent\textbf{Funding: }This study received support from the European Space Agency Climate Change Initiative (ESA-CCI) Tipping Elements SIRENE project (contract no. 4000146954/24/I-LR). N.B. also acknowledges funding by the Volkswagen Foundation. This is ClimTip contribution \#144; the ClimTip project has received funding from the European Union's Horizon Europe research and innovation programme under grant agreement No. 101137601.

\vspace{1em}
\noindent\textbf{Data Availability: }Visit this associated \href{https://github.com/fsuerhoff-art/Statistical-indicators-for-abrupt-transitions-of-dynamical-systems-with-slow-periodic-forcing.git}{GitHub repository} to access the code generating all figures in this manuscript.


\bibliographystyle{unsrt}
\bibliography{article_references}

\end{document}